\documentclass[reqno]{amsart}

\usepackage{amssymb, amsmath, amsthm, mathrsfs, mathtools}

\newtheorem{theorem}{Theorem}[section]
\newtheorem{proposition}[theorem]{Proposition}

\newtheorem{corollary}[theorem]{Corollary}
\newtheorem{fact}[theorem]{Fact}

\theoremstyle{definition}
\newtheorem{definition}[theorem]{Definition}

\theoremstyle{remark}

\newtheorem{remark}[theorem]{Remark}

\numberwithin{equation}{section}

\begin{document}


\title{Intersection Numbers of Geodesic Arcs}


\author{Yoe Alexander Herrera Jaramillo}
\address{Departmento de Ciencias B\'asicas, Coorporaci\'on Universitaria de Investigaci\'on y Desarrollo - UDI,
Bucaramanga, Santander, Colombia}
\email{yoeherrera@gmail.com}
\urladdr{https://sites.google.com/site/yoeherrera/} 




\begin{abstract} For a compact surface $S$ with constant curvature $-\kappa$ (for some $\kappa>0$) and genus $g\geq2$, we show that the tails of the distribution of the normalized intersection numbers $i(\alpha,\beta)/l(\alpha)l(\beta)$ (where $i(\alpha,\beta)$ is the intersection number of the closed geodesics $\alpha$ and $\beta$ and $l(\cdot)$ denotes the geometric length) are estimated by a decreasing exponential function. As a consequence, we find the asymptotic average of the normalized intersection numbers of pairs of closed geodesics on $S$. In addition, we prove that the size of the sets of geodesic arcs whose $T$-self-intersection number is not close to $\kappa T^2/(2\pi^2(g-1))$ is also estimated by a decreasing  exponential function. And, as a corollary of the latter, we obtain a result of Lalley which states that most of the closed geodesics $\alpha$ on $S$ with $l(\alpha)\leq T$ have roughly $\kappa l(\alpha)^2/(2\pi^2(g-1))$ self-intersections, when $T$ is large.

{\bf Classification:} 37d40
\end{abstract}

\maketitle


\section{Introduction}

\paragraph{} Let $S$ be a compact surface of constant curvature $-\kappa$, for some $\kappa>0$, and genus $g\geq2$. A {\it geodesic} ({\it parametrized by the arc length}) {\it on} $S$ is a smooth locally distance-minimizing curve $\gamma:\mathbb R\to S$. For every $x\in S$ and every unit vector $v$ tangent to $S$ at $x$, there is a unique geodesic $\gamma_{(x,v)}$ on $S$ such that $\gamma_{(x,v)}(0)=x$ and $\dot\gamma_{(x,v)}(0)=v$, where $\dot\gamma(t)$ denotes the unit vector tangent to $\gamma$ at $\gamma(t)$. The restriction $\gamma'=\gamma_{\mid[a,b]}$ for $-\infty\leq a<b\leq\infty$ is called a {\it geodesic arc} or {\it segment} and its {\it length} is $l(\gamma')=b-a$. The geodesic $\gamma$ is {\it closed} if there exists $l>0$ such that $\gamma([0,l])=\gamma(\mathbb R)$, and in this case, we say that $l(\gamma)=\min\{l\mid\gamma([0,l])=\gamma(\mathbb R)\}$.

\paragraph{} Two geodesics $\gamma$ and $\eta$ on $S$ are {\it identical} if they both have the same trace, that is, there is $r\neq 0$ such $\gamma(t)=\eta(t+r)$ and $\dot\gamma(t)=\dot\eta(t+r)$, for every $t\in\mathbb R$. Let $[\gamma]$ be the equivalence class formed by all geodesics on $S$ that are identical to $\gamma$. We choose a representative geodesic from each class and form a set that we denote by $\mathbb G$. Let $C\mathbb G$ be the subset of $\mathbb G$ consisting of the geodesics that are closed. Let $C\mathbb G_T=\{\gamma\in C\mathbb G:l(\gamma)\leq T\}$ and $N(T)$ be the cardinality of $C\mathbb G_T$. H. Huber
 proved in~\cite[Theorem 10]{Huber1959:Bibliography} that the number $N(T)$ satisfies the asymptotic formula $N(T)\sim e^{\sqrt\kappa T}/{\sqrt\kappa T}$, that is,$\underset{T\to\infty}\lim\frac{T\sqrt\kappa N(T)}{e^{T\sqrt\kappa}}=1$.

\begin{definition}\label{intnum} Let $T>0$, and $\gamma$ and $\eta$ be geodesics on $S$. The $T${\it-intersection number of} $\gamma$ {\it and} $\eta$ is denoted by $i^T(\gamma,\eta)$ and defined by
$$i^T(\gamma,\eta)=\#\{x\mid\gamma(r)=\eta(t)=x;\dot\gamma(r),\dot\eta(t)\text{ are non-parallel, for some }r,t\in[0,T]\}.$$In particular, $i^T(\gamma,\gamma)$ is the $T${\it-self-intersection number of} $\gamma$. \end{definition}

\begin{remark}If $\gamma$ is a closed geodesic or a geodesic arc with $l(\gamma)\leq T$, then $$i^T(\gamma,\eta)=i(\gamma,\eta),$$ where $i(\gamma,\eta)$ is the ({\it geometric}) {\it intersection number} of $\gamma$ and $\eta$. And, $i(\gamma,\gamma)$ is the {\it self-intersection number of} $\gamma$.
\end{remark}

\paragraph{} The intersection numbers have been of interest to many researchers and here are some of the most relevant results so far achieved. Lalley showed in ~\cite{Lalley1996:Bibliography} that for $T$ large enough, the self-intersection number of most of the closed geodesics $\alpha$ with $l(\alpha)\leq T$ is about $\kappa l(\alpha)^2/(2\pi^2(g-1))$. Later, Pollicott and Sharp generalized this result to self-intersections of closed geodesics with and angle in a given interval (see~\cite{Pollicott2006:Bibliography}). Recently, Chas and Lalley in~\cite{Lalley2011:Bibliography} proved that if a free homotopy class of curves on a surface with boundary is chosen at random from among all classes of word length $m$, then the distribution of the self-intersection numbers appropriately scaled approaches the Gaussian distribution, for $m$ ``large enough." Furthermore, Lalley also showed in~\cite{Lalley2013:Bibliography}
that the random variable $(N_T-\kappa T^2/(2\pi^2(g-1))/T$ has a limit distribution as $T\to\infty$, where $N_T$ is the number of self-intersections of a closed geodesic on $S$ of length $\leq T$ randomly chosen.

\paragraph{} In this paper, we prove that the tails of the distribution of the {\it normalized intersection numbers} of the pairs of elements of $C\mathbb G$, that is $i(\alpha,\beta)/l(\alpha)l(\beta)$ for $\alpha,\beta\in C\mathbb G$, are estimated by a decreasing exponential function.

\begin{theorem}\label{mainthm} Let $\epsilon>0$. There exists $\delta>0$ such that
$$\frac{1}{N(R)N(T)}\#\bigg\{(\alpha,\beta)\in C\mathbb G_R\times C\mathbb G_T:\bigg|\frac{i(\alpha,\beta)}{l(\alpha)l(\beta)}-\frac{\kappa}{2\pi^2(g-1)}\bigg|\geq\epsilon\bigg\}=O\big(e^{-\delta R}\big),$$
as $R\to\infty$, with $T\geq R$. \end{theorem}

\paragraph{} Theorem~\ref{mainthm} allows us to show that the average of the normalized intersection numbers of pairs of closed geodesics of length at most $R$ and $T$ is asymptotically equal to $\kappa/(2\pi^2(g-1))$.

\begin{corollary}\label{average} $\displaystyle{
\frac{1}{N(R)N(T)}\sum_{(\alpha,\beta)\in C\mathbb G_R\times C\mathbb G_T}\frac{i(\alpha,\beta)}{l(\alpha)l(\beta)}\sim \frac\kappa{2\pi^2(g-1)}\text{, as }R,T\to\infty.}$\end{corollary}

\paragraph{} In order to introduce our other results we need the following definitions.

\paragraph{} Let $T^1(S)=\{\mathbf v=(x,v)\mid x\in S,v\in T_x(S),\|v\|=1\}$ be the {\it unit tangent bundle of} $S$. Let $\vartheta$ denote the {\it Riemannian measure on} $T^1(S)$, i.e., the volume measure of $T^1(S)$. In addition, let $\overline\vartheta$ denote the normalized Riemannian measure, that is, $\overline\vartheta=\dfrac1{\vartheta(T^1(S))}\vartheta=\dfrac\kappa{2\pi^2(g-1)}\vartheta$.

\paragraph{} By identifying the unit tangent bundle of $S$ with the set of geodesics on $S$ we prove that the size (or $\overline\vartheta$-measure) of the subset of $T^1(S)$ consisting of vectors whose corresponding geodesics have the normalized $T$-self-intersection number not close to $\kappa/(2\pi^2(g-1))$ is bounded by a decreasing exponential function.

\begin{theorem}\label{thmgeodcurves} Let $\epsilon>0$. There exists $\delta>0$ such that
$$\overline\vartheta\bigg\{\mathbf v\in T^1(S):\bigg|\frac{i^T(\gamma_\mathbf v,\gamma_\mathbf v)}{T^2}-\frac\kappa{2\pi^2(g-1)}\bigg|\geq\epsilon\bigg\}=O\big(e^{-\delta T}\big)\text{, as }T\to\infty.$$
\end{theorem}

\paragraph{} As a consequence of Theorem ~\ref{thmgeodcurves}, we obtain the result given by Lalley in~\cite[Theorem 1]{Lalley1996:Bibliography}.

\begin{corollary}[Lalley]\label{lalley}
For every $\epsilon>0$,
$$\underset{T\to\infty}\lim\frac{1}{N(T)}\#\bigg\{\gamma\in C\mathbb G_T:\bigg|i(\gamma,\gamma)-\frac{\kappa l(\gamma)^2}{2\pi^2(g-1)}\bigg|<\epsilon l(\gamma)^2\bigg\}=1.$$
\end{corollary}

\paragraph{} The outline of this paper is the following. Section~\ref{prelim} is the collection of definitions
and results needed in the demonstrations of Theorems~\ref{mainthm} and~\ref{thmgeodcurves}, and Section~\ref{results}
contains the proofs of these theorems as well as the proofs of Corollary~\ref{average} and
Corollary~\ref{lalley}. For detailed explanation of all the concepts (or a different approach on them)
used in this work, please see~\cite{Bonahon1988:Bibliography},~\cite{Katok1995:Bibliography},~\cite{Kifer1994:Bibliography},~\cite{Parkkonen2013:Bibliography} and~\cite{Pollicott2006:Bibliography}.

{\bf Acknowledgments.} The author of this paper wants to thank both his PhD adviser Dr. Kasra Rafi for all his support and the referee of this paper for his or her helpful suggestions.

\section{Preliminaries}\label{prelim}

\subsection{Measure of Maximum Entropy}

\paragraph{} The map $\varphi:T^1(S)\times\mathbb R\to T^1(S)$ defined by $\varphi(\mathbf v,t)=\varphi^t\mathbf v=(\gamma_{\mathbf v}(t),\dot\gamma_{\mathbf v}(t))$ is the {\it geodesic flow over} $S$. Let $h_{top}(\varphi)$ be the {\it topological entropy of} $\varphi$. A measure $\mu$ on $T^1(S)$ is $\varphi${\it -invariant} if $\mu(\varphi^t(E))=\mu(E)$, for every $t\in\mathbb R$ and every Borel set $E$ of $T^1(S)$. For instance, the measure $\vartheta$ is $\varphi$-invariant. Denote by $\mathscr P_\varphi$ the set of $\varphi$-invariant  probability measures on $T^1(S)$ equipped with the weak*-topology, and for $\mu\in\mathscr P_\varphi$, let $h_\mu(\varphi)$ denote its {\it measure theoretic entropy with respect to} $\varphi$ (please see~\cite[\S 4.3]{Katok1995:Bibliography} for definitions.) The Variational Principal (proven by T.N.T. Goodman) in~\cite{Goodman1971:Bibliography} states that $h_{top}(\varphi)=\sup_{\mu\in\mathscr P_\varphi} h_\mu(\varphi)$. In fact, Bowen proved in~\cite{Bowen1977:Bibliography} that in our case this supremum is actually a maximum and is uniquely achieved by the normalized Riemannian measure $\overline{\vartheta}$, with $h_{\overline\vartheta}(\varphi)=h_{top}(\varphi)=\sqrt\kappa$. Therefore, $\overline\vartheta$ coincides with the {\it measure of maximum entropy on} $T^1(S)$. Moreover, $\overline\vartheta$ also coincides with the {\it Margulis-Bowen measure} from ~\cite{Margulis2004:Bibliography}, and in this work, we use the characterization of this measure given by Bowen in~\cite{Bowen1972:Bibliography}.

\paragraph{} The ($\varphi$-){\it orbit of} $\mathbf v\in T^1(S)$ is the set $\{\varphi^t\mathbf v\mid t\in\mathbb R\}$. These orbits form a partition of $T^1(S)$. Note that there is a one-to-one correspondence between the set of orbits and the set $\mathbb G$. The vector $\mathbf v\in T^1(S)$ and its orbit are {\it periodic} if there exists $l>0$ such that $\varphi^l\mathbf v=\mathbf v$, the number $l$ is a {\it period} and the minimal period is precisely $l(\gamma_\mathbf v)$.

\paragraph{} For a periodic orbit $\gamma$, Bowen defined the {\it occupation measure} $\zeta_\gamma$ on $T^1(S)$ by
\begin{eqnarray}\label{bowen2}\zeta_\gamma(E)=\int^{l(\gamma)}_0\chi_E(\varphi^t\mathbf v)dt,\end{eqnarray}
for $\mathbf v\in\gamma$ and $E$ a Borel set of $T^1(S)$. In addition, Bowen proved in ~\cite[(5.5)]{Bowen1972:Bibliography} the following.

\begin{theorem}[Bowen]\label{bowen}
The periodic orbits of the geodesic flow $\varphi$ are equidistributed with respect to the measure of maximium entropy $\overline\vartheta$ as the period tends to $+\infty$. More precisely, for any Borel set $E$ with $\overline\vartheta(\partial E)=0$,
$$\overline\vartheta(E)=\underset{T\to\infty}\lim\frac1{N(T)}\sum_{\gamma\in C\mathbb G_T}\frac{\zeta_\gamma}{l(\gamma)}(E).$$
\end{theorem}

\subsection{Geodesic Currents}

\paragraph{} Let $\mathbb L=T^1(S)/\sim$, with $(x,v)\sim(x,-v)$, be the {\it line bundle of} $S$ and $\mathcal F$ be the {\it foliation} of $\mathbb L$ by $\varphi$-orbits.
A ({\it geodesic}) {\it current} $\mu$ {\it on} $S$ is a positive transverse invariant
measure for the geodesic foliation $\mathcal F$. The set of currents on $S$ equipped with the weak*topology is denoted by $\mathcal C$ and called the {\it space of currents on} $S$.

\paragraph{} Given any $\varphi$-invariant measure $\mu$, we can consider the associated
transverse measure $\widetilde\mu$ for the foliation $\mathcal F$. Each $\widetilde\mu\in\mathcal C$ is normalized
by the requirement that (locally) $\mu=\widetilde\mu\times dt$, where $dt$ is the one--dimensional Lebesgue measure along leaves in $\mathcal F$. The current associated to $\vartheta$ is called the {\it Liouville current on} $S$. In this paper, we identify the measure $\mu$ with the current $\tilde{\mu}$.

\paragraph{} The basic example of a current is the one associated to  a closed geodesic $\gamma$ on $S$. To this geodesic $\gamma$ corresponds a compact leaf $\widetilde\gamma$ of $\mathcal F$. We associate to it
the current $\mu_\gamma$ which induces on each transverse manifold $V$ the Dirac
measure at the point $V\cap\widetilde\gamma$. Such current corresponds to the $\varphi$-invariant measure $\zeta_\gamma$, as defined in (\ref{bowen2}).

\paragraph{} Observe that it is always possible to add two geodesic currents, and to multiply a geodesic current by a non-negative real number. Then, the space $\mathcal C$ appears as the completion of the space of real multiples of homotopy
classes of closed curves by the following fact, which we state although we will
not make use of it in this paper.

\begin{proposition} The uniform space $\mathcal C$ is complete, and the real multiples of homotopy classes of closed curves are dense in it. \end{proposition}

\subsection{Intersection Form}

\paragraph{} Starting from the bundle $\mathbb L\to S$, we can consider the Whitney sum
$\mathbb L\oplus\mathbb L\to S$. In other words, $\mathbb L\oplus\mathbb L$ is the 4-–dimensional manifold
of triples $(x,\lambda_1, \lambda_2)$, where $x\in S$ and $\lambda_1$ and $\lambda_2$ are two lines in the tangent
space $T_x(S)$. Forgetting the first or the second line defines two projections
$p_1$ and $p_2$ from $\mathbb L\oplus\mathbb L$ to $\mathbb L$. We consider the two foliations $\mathcal F_1$ and
$\mathcal F_2$ of codimension $2$ in $\mathbb L\oplus\mathbb L$, whose leaves are the preimages of the
leaves of $\mathcal F$ by, respectively, $p_1$ and $p_2$. These foliations
are transverse outside the diagonal $\triangle=\{(x,v,v)\mid (x,v)\in\mathbb L\}$ of $\mathbb L\oplus\mathbb L$.

\paragraph{} Let $\mu$ and $\nu$ be two currents. Through $p_1$, $\mu$ induces a transverse
invariant measure $\widehat\mu_1$
on $\mathcal F_1$, which, by transversality of $\mathcal F_1$ and $\mathcal F_2$,
gives outside $\triangle$ a measure on each leaf of $\mathcal F_2$. Similarly, $\nu$ induces outside
$\triangle$ a measure $\widehat\nu_2$
on each leaf of $\mathcal F_1$. Consider then the product measure $\widehat\mu_1\times\widehat\nu_2$ on $\mathbb L\oplus\mathbb L\setminus\triangle$. The total mass of this measure is finite. The {\it intersection form of} $\mu$ {\it and} $\nu$ is $$\imath(\mu,\nu)=\widehat\mu_1\times\widehat\nu_2(\mathbb L\oplus\mathbb L\setminus\triangle).$$

\begin{remark}The {\it normalized Liouville current} denoted by $\overline\vartheta$ (which corresponds to the normalized Riemannian measure) is defined by $$\imath(\vartheta,\overline\vartheta)=1.$$
\end{remark}

\paragraph{} By identifying the closed geodesic $\alpha$ on $S$ with the current $\zeta_\alpha$, Bonahon proved the following facts in~\cite[Theorem 4.1]{Bonahon1986:Bibliography} and~\cite[Proposition 15]{Bonahon1988:Bibliography}.

\begin{theorem}[Bonahon]\label{bonahon} The intersection form function $\imath:\mathcal C\times\mathcal C\to[0,\infty)$ is a continuous extension of the intersection number function. In particular, for $\alpha$ and $\beta$ closed geodesics on $S$, $$\imath(\zeta_\alpha,\zeta_\beta)=i(\alpha,\beta).$$
\paragraph{} In addition, $$\imath(\vartheta,\vartheta)=\dfrac{2\pi^2(g-1)}\kappa=\frac{\vartheta(T^1(S))}2$$
and
$$\imath(\overline\vartheta,\overline\vartheta)=\dfrac\kappa{2\pi^2(g-1)}=\frac2{\vartheta(T^1(S))}.$$
\end{theorem}

\section{Results}\label{results}

\paragraph{} The proofs of Theorem~\ref{mainthm} and~\ref{thmgeodcurves} are based on both the continuity of the intersection form function (which is continuous by Theorem~\ref{bonahon}) and a deviation result given by Y. Kifer ~\cite{Kifer1994:Bibliography}.

\paragraph{} For $T>0$ and $\mathbf v\in T^1(S)$, Y. Kifer defined the occupation measure $\zeta^T_\mathbf v$ by
$$\zeta^T_{\mathbf v}(E)=\int^T_0\chi_E(\varphi^t\mathbf v)dt,$$
for every Borel set $E$ of $T^1(S)$.

\paragraph{} Note that if $\gamma$ is a periodic orbit, we have $\zeta^{l(\gamma)}_{\mathbf v}=\zeta_\gamma$, for $\mathbf v\in\gamma$.

\begin{fact}\label{factg} The intersection form can be extended to the whole set of (positive) finite measures (not necessarily $\varphi$-invariant). By abuse of notation, we denote this extension also by $\imath$. Such extension satisfies the following:
for $\mathbf v,\mathbf w\in T^1(S)$ and $T>0$, $$\imath\Big(\zeta^T_\mathbf v,\zeta^T_\mathbf w\Big)=i^T(\gamma_\mathbf v,\gamma_\mathbf w).$$ \end{fact}

\paragraph{} Since $T^1(S)$ is a compact metric space and $\varphi$ is a hyperbolic dynamical system, the deviation results of Y. Kifer in~\cite[Theorem 3.4]{Kifer1990:Bibliography} and ~\cite[Theorem 2.1]{Kifer1994:Bibliography}, respectively, can be translated into our setting in the following way.
\begin{theorem}[Kifer]\label{kifer} For any closed subset $K$ of $\mathscr P$, the space of probability measures on $T^1(S)$,
\begin{eqnarray*}\underset{T\to\infty}\limsup\frac{1}{T}
\log\overline\vartheta\bigg\{\mathbf v\in T^1(S):\frac{\zeta^T_\mathbf v}T\in K\bigg\}\leq-\underset{\mu\in K}\inf I(\mu), \end{eqnarray*}
where
\[I(\mu)=
\begin{dcases} h_{top}(\varphi)-h_\mu(\varphi), & \mu\in\mathscr P_\varphi\\
\infty, & otherwise
\end{dcases}.\]
\end{theorem}

\begin{theorem}[Kifer]\label{kif} Let $\mathcal U$ be an open neighborhood of the measure of maximal entropy $\overline\vartheta$ in the set of $\varphi$-invariant probability measures on $T^1(S)$. Then
\begin{eqnarray*}\underset{T\to\infty}\lim\frac{1}{N(T)}
\#\{\gamma\in C\mathbb G_T:\zeta_\gamma/l(\gamma)\not\in\mathcal U\}=O\big(e^{-\delta T}\big), \end{eqnarray*}
as $T\to\infty$, where
$\delta=\underset{\mu\in\mathcal U^c}\inf\{h_{top}(\varphi)-h_\mu(\varphi)\}$.
\end{theorem}

\begin{proof}[Proof of Theorem~\ref{thmgeodcurves}] Let $\epsilon>0$. Consider the set
\begin{eqnarray}\label{setk}K:=\{\mu\in\mathscr P:|\imath(\mu,\mu)-\imath(\overline\vartheta,\overline\vartheta)|\geq\epsilon\}.\end{eqnarray}

\paragraph{} By Theorem~\ref{bonahon} and Fact~\ref{factg},
\begin{eqnarray}\label{geodcurves} \bigg\{\mathbf v\in T^1(S):\bigg|\frac{i^T(\gamma_\mathbf v,\gamma_\mathbf v)}{T^2}-\frac\kappa{2\pi^2(g-1)}\bigg|\geq\epsilon\bigg\}=\bigg\{\mathbf v\in T^1(S):\frac{\zeta^T_\mathbf v}T\in K\bigg\}.\end{eqnarray}

\paragraph{} By (\ref{geodcurves}), it is enough to prove that there exists $\delta>0$ such that $$\overline\vartheta\bigg\{\mathbf v\in T^1(S):\frac{\zeta^T_\mathbf v}T\in K\bigg\}=O\big(e^{-\delta T}\big),\text{ as }T\to\infty.$$

\paragraph{} The intersection form function $\imath$ is continuous by Theorem~\ref{bonahon}, then, the set $K$ in (\ref{setk}) is a closed subset of $\mathscr P$. Therefore, by Theorem~\ref{kifer},
$$\underset{T\to\infty}\limsup\frac{1}{T}
\log\overline\vartheta\bigg\{\mathbf v\in T^1(S):\frac{\zeta^T_\mathbf v}T\in K\bigg\}\leq-\underset{\mu\in K}\inf I(\mu).$$

\paragraph{} If $K\cap\mathscr P_\varphi=\emptyset$, then  $\underset{\mu\in K}\inf I(\mu)=\infty$. Thus, $$\underset{T\to\infty}\limsup\frac{1}{T}
\log\overline\vartheta\bigg\{\mathbf v\in T^1(S):\frac{\zeta^T_\mathbf v}T\in K\bigg\}\leq-\infty.$$ Hence $\overline\vartheta\Big\{\mathbf v\in T^1(S):\zeta^T_\mathbf v\in K\Big\}\leq e^{-\infty}=0=O\big(e^{-\delta T}\big),$ as $T\to\infty$, for any $\delta>0$.

\paragraph{} If $K\cap\mathscr P_\varphi\neq\emptyset$, given that $\overline\vartheta$ is the unique probability measure of $T^1(S)$ with maximum entropy $h_{\overline\vartheta}(\varphi)=h_{top}(\varphi)=\sqrt\kappa$, then $$0<\delta=\underset{\mu\in K}\inf I(\mu)=\underset{\mu\in K\cap\mathscr P_\varphi}\inf I(\mu)=\underset{\mu\in K\cap\mathscr P_\varphi}\inf(\sqrt\kappa-h_\mu(\varphi)) <\infty$$
is such that
$$\underset{T\to\infty}\limsup\frac{1}{T}
\log\overline\vartheta\bigg\{\mathbf v\in T^1(S):\frac{\zeta^T_\mathbf v}T\in K\bigg\}\leq-\delta.$$
\paragraph{}Hence,$$\overline\vartheta\bigg\{\mathbf v\in T^1(S):\frac{\zeta^T_\mathbf v}T\in K\bigg\}=O\big(e^{-\delta T}\big),\text{ as }T\to\infty.$$
\end{proof}
\begin{proof}[Proof of Corollary~\ref{lalley}] Let $T,\epsilon>0$ and 
$$\mathscr O(T,\epsilon):=\bigg\{\gamma\in C\mathbb G_T:\bigg|\frac{i(\gamma,\gamma)}{l(\gamma)^2}-\frac\kappa{2\pi^2(g-1)}\bigg|<\epsilon\bigg\}.$$

\paragraph{} Consider $$\mathcal U=K^c=\{\mu\in\mathscr P:|\imath(\mu,\mu)-\imath(\overline\vartheta,\overline\vartheta)|<\epsilon\}.$$

\paragraph{} Then, $\mathscr O(T,\epsilon)=\{\gamma\in C\mathbb G_T:\zeta_\gamma\in\mathcal U\}.$

\paragraph{} By Theorem 3.3, we have $$\frac{1}{N(T)}\#\{\gamma\in C\mathbb G_T:\zeta_\gamma/l(\gamma)\not\in\mathcal U\}=O\big(e^{-\delta T}\big).$$

\paragraph{} Consequently,
$$\underset{T\to\infty}\lim\frac{1}{N(T)}\#\bigg\{\gamma\in C\mathbb G_T:\bigg|i(\gamma,\gamma)-\frac{\kappa l(\gamma)^2}{2\pi^2(g-1)}\bigg|<\epsilon l(\gamma)^2\bigg\}=1.$$
\end{proof}
\begin{proof}[Proof of Theorem~\ref{mainthm}] Consider the function $\imath:\mathscr P\times\mathscr P\to\imath(\mathscr P\times\mathscr P)$. This function is continuous since it is the restriction of the intersection form function $\imath$, which continuous by Theorem~\ref{bonahon}, to $\mathscr P\times\mathscr P$ a closed subset of $\mathcal  C\times\mathcal C$.

\paragraph{} Therefore, for  $\epsilon>0$, the set $\mathcal Z=\mathfrak\imath^{-1}\bigg(\dfrac{\kappa}{2\pi^2(g-1)}-\epsilon,\dfrac{\kappa}{2\pi^2(g-1)}+\epsilon\bigg)$ is an open subset of $\mathscr P\times\mathscr P$ since it is the preimage under $\imath$ of the ball of radius $\epsilon$ centered at $\imath(\overline\vartheta,\overline\vartheta)=\kappa/(2\pi^2(g-1))$.

\paragraph{} Let $R,T>0$ with $R\leq T$ and 
\begin{eqnarray*}\label{setw}\mathcal W_{R,T}=\bigg\{(\alpha,\beta)\in C\mathbb G_R\times C\mathbb G_T:\bigg(\frac{\zeta_\alpha}{l(\alpha)},\frac{\zeta_\beta}{l(\beta)}\bigg)\in\mathcal Z\bigg\}.\end{eqnarray*}

\paragraph{} By Theorem~\ref{bonahon},
\begin{eqnarray}\label{thmsub1}\bigg\{(\alpha,\beta)\in C\mathbb G_R\times C\mathbb G_T:\bigg|\frac{i(\alpha,\beta)}{l(\alpha)l(\beta)}-\frac{\kappa}{2\pi^2(g-1)}\bigg|\geq\epsilon\bigg\}=C\mathbb G_R\times C\mathbb G_T\setminus\mathcal W_{R,T}.
\end{eqnarray}

\paragraph{} Hence, by (\ref{thmsub1}), it is enough to prove that there exists $\delta>0$ such that
$$\frac{\#C\mathbb G_R\times C\mathbb G_T\setminus \mathcal W_{R,T}}{N(R)N(T)}=O\big(e^{-\delta R}\big).$$

\paragraph{} Since $(\overline\vartheta,\overline\vartheta)\in\mathcal Z$ and $\mathcal Z$ is an open set of the product topology of $\mathscr P\times\mathscr P$, there exist $\mathcal U,\mathcal V\subseteq\mathscr P$ open neighborhoods of $\overline\vartheta$ in $\mathscr P$ such that $\mathcal U\times\mathcal V\subseteq\mathcal Z$.

\paragraph{} Let $\mathcal U_R=\bigg\{\alpha\in C\mathbb G_R:\dfrac{\zeta_\alpha}{l(\alpha)}\in\mathcal U\bigg\}$ and $\mathcal V_T=\bigg\{\beta\in C\mathbb G_T:\dfrac{\zeta_\beta}{l(\beta)}\in\mathcal V\bigg\}$. 

\paragraph{} Given that both $\mathcal U$ and $\mathcal V$ are open neighborhoods of $\overline\vartheta$ on $\mathscr P$, Theorem~\ref{kifer}, guarantees the existence of $\delta_1,\delta_2>0$ depending on $\mathcal U$ and $\mathcal V$, respectively, such that
$$\frac{\#C\mathbb G_R\setminus\mathcal U_R}{N(R)}=\frac1{N(R)}\#\{\gamma\in C\mathbb G_R:\zeta_\gamma/l(\gamma)\not\in\mathcal U\}=O\big(e^{-\delta_1R}\big)$$ and $$\frac{\#C\mathbb G_T\setminus\mathcal V_T}{N(T)}=\frac1{N(T)}\#\{\gamma\in C\mathbb G_T:\zeta_\gamma/l(\gamma)\not\in\mathcal V\}=O\big(e^{-\delta_1T}\big).$$

\paragraph{} Thus, since $\mathcal U_R \times\mathcal V_T \subseteq \mathcal W_{R,T}$, we get as $R,T\to\infty$,
\begin{align*}\frac{\#C\mathbb G_R\times C\mathbb G_T\setminus \mathcal W_{R,T}}{N(R)N(T)}&\leq\frac{\#C\mathbb G_R\times C\mathbb G_T\setminus\mathcal U_R\times\mathcal V_T}{N(R)N(T)}\\
&\leq\frac{\#C\mathbb G_R\setminus \mathcal U_R\cdot\#C\mathbb G_T\setminus\mathcal V_T}{N(R)N(T)}\\
&\quad\quad\quad +\frac{\# C\mathbb G_R \setminus\mathcal U_R\cdot\#\mathcal V_T}{N(R)N(T)}+\frac{\#C\mathbb G_T\setminus\mathcal V_T\cdot\#\mathcal U_R}{N(R)N(T)} \\ &=O\big(e^{-\delta_1R}\big)O\big(e^{-\delta_2T}\big)+O\big(e^{-\delta_1R}\big)+O\big(e^{-\delta_2T}\big)=O\big(e^{-\delta_1R}\big). \end{align*}
\end{proof}

\paragraph{} For the proof of Corollary~\ref{average}, we need a bound for the intersection number of pairs of closed geodesics on $S$. Here, we provide a universal bound for the normalized intersection numbers of pairs of closed geodesics. It is worth noting that such bound can also be deduced by the techniques used by A. Basmajian in~\cite{Basmajian2010:Bibliography}.

\paragraph{} The {\it injectivity radius at} a point $x\in S$ is the largest radius for which the exponential map at $x$ is a diffeomorphism. The {\it injectivity radius of} $S$, which we denote by $\varrho$, is the infimum of the injectivity radii of all points of $S$. By the definition of $\varrho$, the least length of an essential loop on $S$ is $2\varrho$.

\begin{proposition}\label{bound} Let $\alpha$ and $\beta$ be two closed geodesics on $S$.
Then $$\frac{i(\alpha,\beta)}{l(\alpha)l(\beta)}\leq\frac1{\varrho^2}.$$
\end{proposition}

\begin{proof} Let $\bar\alpha$ be a sub-arc of $\alpha$ with $l(\bar\alpha)<\varrho$ and such that $i(\bar\alpha,\beta)\geq i(\alpha^*,\beta)$ for any sub-arc $\alpha^*$ of $\alpha$ with $l(\alpha^*)<\varrho$.
Hence, \begin{eqnarray}\label{boundint}i(\alpha,\beta)\leq\bigg\lceil\frac{l(\alpha)}\varrho\bigg\rceil i(\bar\alpha,\beta).\end{eqnarray}

\paragraph{} Let $\{x_1,\ldots,x_n\}$ be the ordered set of points of intersection of $\bar\alpha$ and $\beta$, with $\beta^{-1}(x_i)\leq\beta^{-1}(x_{i+1})$, for $1\leq i\leq n-1$,
and $n=i(\bar\alpha,\beta)$. Let $\beta_k$ be the sub-arc of $\beta$ from $x_k$ to $x_{k+1}$, for $1\leq k\leq n-1$,
and, $\beta_n$ be the sub-arc of $\beta$ from $x_n$
to $x_1$. Similarly,
let $\bar\alpha_k$ be the sub-arc of $\bar\alpha$ from $x_k$
to $x_{k+1}$, for $1\leq k\leq n-1$, and $\bar\alpha_n$
be the sub-arc of $\bar\alpha$ from $x_n$ to $x_1$. 

\paragraph{} Consider $\gamma_k$ the concatenation of $\bar\alpha_k$ and $\beta_k$, for $1\leq k\leq n$. Thus, $\gamma_k$ is an essential loop of $S$, for $1\leq k \leq n$.

\paragraph{} Hence, $2\varrho\leq l(\gamma_k)=l(\bar\alpha_k)+l(\beta_k)\leq l(\bar\alpha)+l(\beta_k)<\varrho+l(\beta_k)$, which implies $\varrho<l(\beta_k)$, for $1\leq k\leq n$.

\paragraph{} Consequently, $n\varrho=n\sum^n_{i=1}\varrho<\sum^n_{i=1}l(\beta_k)\leq l(\beta)$. Therefore, $n=i(\bar\alpha,\beta)<\dfrac{l(\beta)}{\varrho}$. Thus, by (\ref{boundint}), we conclude 
$\displaystyle{
(\alpha,\beta)\leq \left\lceil\frac{l(\alpha)}{\varrho}\right\rceil i(\bar \alpha,\beta)\leq\frac{l(\alpha)}{\varrho} \frac{l(\beta)}{\varrho}= \frac{l(\alpha)l(\beta)}{\varrho^2}.}$ \end{proof}

\begin{proof}[Proof of Corollary~\ref{average}] Let $\epsilon>0$. For $R,T>0$ with $R\leq T$, consider the set $\mathcal W_{R,T}$ defined in (\ref{setw}) from the proof of Theorem~\ref{mainthm}. In addition, let $\delta,J,C>0$ be constants satisfying the conclusion of such theorem, that is, for $J\leq R\leq T$, we have
$$\frac{\#C\mathbb G_R\times C\mathbb G_T\setminus\mathcal W_{R,T}}{N(R)N(T)}\leq\frac{C}{e^{\delta R}}.$$

\paragraph{} Moreover, let $J$ be such that $Ce^{-\delta R}<\epsilon$, whenever $R>J$.

\paragraph{} By Proposition~\ref{bound}, we have $\underset{C\mathbb G_R\times C\mathbb G_T}\sup\dfrac{i(\alpha,\beta)}{l(\alpha)l(\beta)}\leq\dfrac1{\varrho^2}$. Therefore, for $J<R\leq T$, we have
\begin{align*}\Bigg|\frac{2\pi^2(g-1)}{\kappa N(R)N(T)}&\Bigg(\sum_{\substack{\alpha\in C\mathbb G_R \\ \beta\in C\mathbb G_T}}\frac{i(\alpha,\beta)}{l(\alpha)l(\beta)}\Bigg)-1\Bigg|\\
\quad\quad\quad\quad\quad\quad&=\Bigg|\frac{2\pi^2(g-1)}{\kappa N(R)N(T)}\sum_{\substack{\alpha\in C\mathbb G_R \\ \beta\in C\mathbb G_T}}\Bigg(\frac{i(\alpha,\beta)}{l(\alpha)l(\beta)}-\frac{\kappa}{2\pi^2(g-1)}\Bigg)\Bigg|\\
&\leq\frac{2\pi^2(g-1)}{\kappa N(R)N(T)}\Bigg(\sum_{(\alpha,\beta)\in\mathcal W_{R,T}}\Bigg|\frac{i(\alpha,\beta)}{l(\alpha)l(\beta)}-\frac{\kappa}{2\pi^2(g-1)}\Bigg|\\
&\quad\quad\quad\quad\quad+\sum_{(\alpha,\beta)\in C\mathbb G_R\times C\mathbb G_T\setminus\mathcal W_{R,T}}\Bigg|\frac{i(\alpha,\beta)}{l(\alpha)l(\beta)}-\frac{\kappa}{2\pi^2(g-1)}\Bigg|\Bigg)\\
&\leq \frac{2\pi^2(g-1)}{\kappa N(R)N(T)}\Bigg(\#\mathcal W_{R,T}\cdot\epsilon\\
&\quad\quad\quad+\#C\mathbb G_R\times C\mathbb G_T\setminus\mathcal W_{R,T}\cdot\underset{\substack{\alpha\in C\mathbb G_R \\ \beta\in C\mathbb G_T}}\sup\bigg|\frac{i(\alpha,\beta)}{l(\alpha)l(\beta)}-\frac{\kappa}{2\pi^2(g-1)}\Bigg|\Bigg)\\
&<\frac{2\pi^2(g-1)}{\kappa} \bigg(\epsilon+\frac{C}{e^{\delta R}}\bigg[\frac{1}{\varrho^2}+\frac{\kappa}{2\pi^2(g-1)}\bigg]\bigg)\\ &<\frac{2\pi^2(g-1)}{\kappa}\bigg(1+\frac{1}{\varrho^2}+\frac{\kappa}{2\pi^2(g-1)}\bigg)\epsilon. \end{align*}

\paragraph{} Given that $\epsilon$ was chosen arbitrarily, we conclude that
$$\underset{R,T\to\infty}\lim\frac{2\pi^2(g-1)}{\kappa N(R)N(T)}\sum_{(\alpha,\beta)\in C\mathbb G_R\times C\mathbb G_T}\frac{i(\alpha,\beta)}{l(\alpha)l(\beta)}=1,$$
or equivalently,
$$\frac{1}{N(R)N(T)}\sum_{(\alpha,\beta)\in C\mathbb G_R\times C\mathbb G_T}\frac{i(\alpha,\beta)}{l(\alpha)l(\beta)}\sim \frac\kappa{2\pi^2(g-1)},$$
as $R,T\to\infty$.\end{proof}

\end{document}